\documentclass{amsart}
\usepackage{amsmath,amssymb,bbm,latexsym,amsthm}
\usepackage[all]{xy}
\usepackage{graphicx}

\newtheorem{prop}{Proposition}[section]
\newtheorem{lemma}[prop]{Lemma}
\newtheorem{defin}[prop]{Definition}
\newtheorem{cor}[prop]{Corollary}
\newtheorem{theorem}[prop]{Theorem}
\newtheorem{claim}[prop]{Claim}

\newcounter{claim}
\newcounter{claim2}
\newcounter{claim3}
\newcounter{claim4}
\newcounter{cor1}
\newcounter{cor2}
\newcounter{cor3}
\newcounter{cor4}

\newcommand{\ho}{\mathcal{O}}
\newcommand{\altZ}{\mathcal{Z}}
\newcommand{\altA}{\mathcal{A}}
\newcommand{\altD}{\mathcal{D}}
\newcommand{\altT}{\mathcal{T}}

\newcommand{\altG}{\mathcal{G}}

\newcommand{\altF}{\mathcal{F}}

\newcommand{\altI}{\mathcal{I}}

\newcommand{\altP}{\mathcal{P}}

\newcommand{\Z}[4]{\left( \begin{array}{cc} #1 & #2 \\ #3 & #4 \end{array} \right)}
\newcommand{\Vek}[2]{\left( \begin{array}{c} #1  \\ #2  \end{array} \right)}
\newcommand{\sur}{\twoheadrightarrow}

\newcommand{\dd}{^{\vee\vee}}

\DeclareMathOperator{\Hom}{Hom}
\DeclareMathOperator{\Ext}{Ext}
\DeclareMathOperator{\stab}{Stab}

\DeclareMathOperator{\Coh}{Coh}
\DeclareMathOperator{\Cohp}{Coh_{(p)}(X)}
\DeclareMathOperator{\Ima}{Im}
\DeclareMathOperator{\Rea}{Re}

\DeclareMathOperator{\rk}{rk}
\DeclareMathOperator{\Id}{Id}
\DeclareMathOperator{\Der}{D}
\DeclareMathOperator{\Supp}{Supp}

\DeclareMathOperator{\Glt}{\widetilde{GL}{}^+(2,\mathbb{R})}
\DeclareMathOperator{\Gl2}{GL^+(2,\mathbb{R})}
\DeclareMathOperator{\Mat2}{Mat(2,\mathbb{R})}
\DeclareMathOperator{\Pic}{Pic}
\DeclareMathOperator{\coker}{coker}

\DeclareMathOperator{\ch}{ch}
\DeclareMathOperator{\ce}{c}
\DeclareMathOperator{\Ho}{H}

\title{Stability conditions on generic complex tori}
\author{Sven Meinhardt}

\begin{document}

\maketitle

\begin{abstract}
In this paper we describe a simply connected component of the complex manifold $\stab(X)$ of stability conditions on a generic complex torus $X$. A generic complex torus is a complex torus $X$ with $\Ho^{p,p}(X)\cap \Ho^{2p}(X,\mathbb{Z})=0$ for all $0<p<\dim X$. 
\end{abstract}

\section{Introduction}

In his paper \cite{Bridgeland02} T. Bridgeland introduced the notion of a stability condition on a triangulated category. His main result states that the space $\stab(X)$ of numerical locally finite stability conditions on the bounded derived category $\Der^b(X)$ of a compact complex manifold $X$ has a natural structure of a complex manifold. During the last years people have tried to describe the manifold $\stab(X)$ for various $X$. The case of curves was treated by Bridgeland \cite{Bridgeland02}, S. Okada \cite{Okada} and E. Macr\`{i} \cite{Macri}. In their paper \cite{HuybStelMac} D. Huybrechts, P. Stellari and E. Macr\`{i} gave a full description for generic K3 surfaces und generic complex tori of dimension two. The condition `generic' means $\Ho^{1,1}(X)\cap\Ho^{2}(X,\mathbb{Z})=0$. T. Bridgeland considered the case of projective K3 surfaces and abelian surfaces in \cite{Bridgeland03}. For these projective surfaces the structure of the space $\stab(X)$ is only partially known. \\
In this paper we construct stability conditions on generic complex tori of any dimension. A complex torus is called generic if
\[ \Ho^{p,p}(X)\cap \Ho^{2p}(X,\mathbb{Z})=0\qquad \forall \;0<p<\dim X.\]
The main result of this paper is the following theorem.
\begin{theorem}
Assume $X$ is a generic complex torus of dimension $d$. Let $U(X)$ be the set of all numerical locally finite stability conditions $\sigma=(Z,\altP)$ such that there exist certain real numbers $\phi$ and $\psi$ such that $k(y)\in \altP(\phi)$ for all $y\in X$ and $L\in \altP(\psi)$ for all $L\in\Pic^0(X)$. Then $U(X)$  is a simply connected component of $\stab(X)$. \\ 
Furthermore, $U(X)$ can be written as a disjoint union of $\Glt$-orbits
\[ U(X)= \bigcup_{0\le p<d}\sigma_{(p)}\cdot \Glt \;\;\;\cup  \bigcup_{1\le p<d \atop \gamma \in (0,1/2)} \sigma_{(p)}^\gamma \cdot\Glt  \]
with explicitly given stability conditions $\sigma_{(p)}$ and $\sigma_{(p)}^\gamma$.
\end{theorem}
Since the case $\dim X=d\le 2$ has already been studied, we restrict ourself to tori of dimension $d \ge 3$. In contrast to the case $d\le 2$ the space $U(X)$ is no longer a covering of its image under the map $\altZ:U(X)\ni \sigma=(Z,\altP) \longmapsto Z\in \altZ(U(X))\subseteq \Ho^\ast(X,\mathbb{C})^\vee$. Furthermore, also in the case $d\ge 3$ it is still open whether or not $\stab(X)$ is connected. \\
Note that the characterizing condition of $U(X)$ is invariant under the Fourier--Mukai transform with respect to the Poincar\'{e} bundle. Hence there is a natural isomorphism $U(X) \cong U(\hat{X})$, where $\hat{X}=\Pic^0(X)$ is the dual torus. \\
The picture below illustrates $U(X)$ and $\altZ(U(X))$ of a generic torus of dimension $d=5$. Note that a point in the helix represents a simply connected 2-dimensional subspace in the $\Glt$-orbit of some stability condition, whereas a point in the circle below represents a 2-dimensional subspace in the $\Gl2$-orbit with the fundamental group $\mathbb{Z}$. 
\begin{center}
\begin{picture}(235,360)
\put(0,0){\includegraphics[width=8cm]{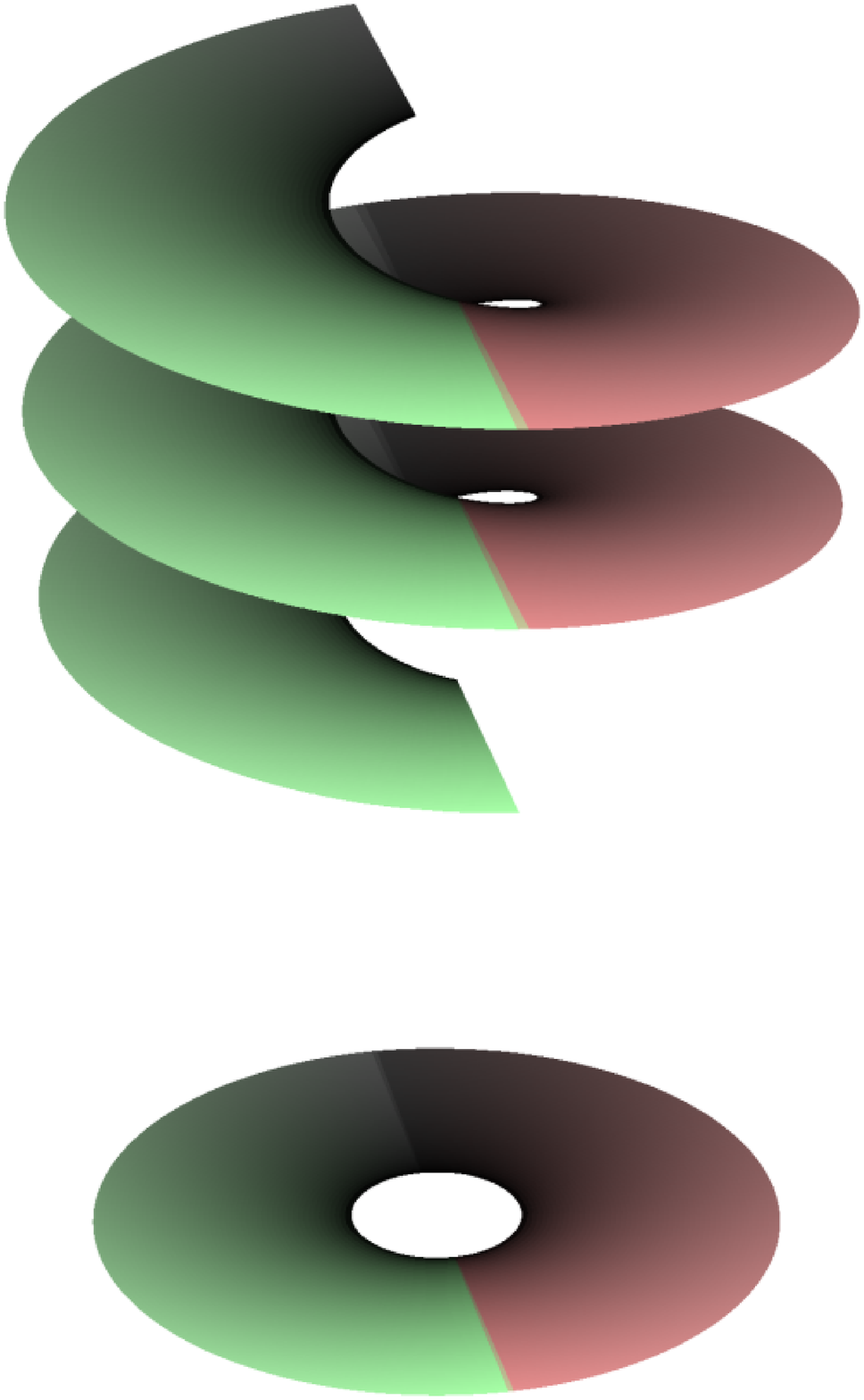}}
\put(213,184){\makebox(0,0){$\cong \Glt$}}
\put(192,95){\makebox(0,0){$\cong \Gl2$}}
\put(125,122){\makebox(0,0){$\altZ$}}
\put(30,144){\makebox(0,0){$\sigma_{(d-1)}$}}
\put(30,320){\makebox(0,0){$\sigma_{(0)}$}}
\put(212,293){\makebox(0,0){$\sigma_{(1)}$}}
\put(234,250){\makebox(0,0){$\sigma_{(1)}^\gamma$}}
\put(40,320){\vector(2,-3){30}}  
\put(200,293){\vector(-2,-1){40}}
\put(222,250){\vector(-4,0){93}}
\put(180,186){\vector(-1,1){20}}
\put(180,182){\vector(-4,1){80}}
\put(33,153){\vector(3,2){36}}
\put(117,133){\vector(0,-1){30}}
\put(160,92){\vector(-1,-1){15}}
\put(160,94){\vector(-2,-1){80}}
\put(73,270){\circle*{0.1}}
\put(157,271){\circle*{0.1}}
\put(125,250){\circle*{0.1}}
\put(73,180){\circle*{0.1}}
\put(73,220){\circle*{0.1}}
\put(157,225){\circle*{0.1}}
\end{picture}
\vspace{1cm}
\end{center}

At this point it still seems very difficult to construct stability conditions on the bounded derived category of projective manifolds of dimension $d\ge 3$. \\
\\
\textbf{Acknowledgements.} I would like to thank Daniel Huybrechts for introducing me to the subject and a lot of fruitful discussions.

\section{Sheaves on generic tori}

In this section we study sheaves on a generic torus $X$ of dimension $d\ge 3$. The following facts and arguments are well known (see e.g. \cite{Verb2} or \cite{Voisin02}). The main result states that on such a torus every reflexive sheaf is locally free and possesses a filtration, whose quotients are line bundles in $\Pic^0(X)$. 
\begin{defin}
A compact complex torus $X$ of dimension $d$ is called generic, if
\[ \Ho^{p,p}(X)\cap \Ho^{2p}(X,\mathbb{Z})=0 \qquad \forall \:0<p<d.\]
\end{defin}
As an  immediate consequence of the definition we get
\begin{itemize}
\item $\Pic(X)=\Pic^0(X)$,
\item the support of any torsion sheaf is a finite set of points in $X$.
\end{itemize}
The last observation leads to the simple but frequently used formula
\begin{equation} \label{eq} \Ext^i(T,F)=\Ext^{d-i}(F,T)^\vee=\Ho^{d-i}(X,T\otimes F^\vee)^\vee=0 \qquad \forall \:i<d
\end{equation}
for a torsion sheaf $T$ and a locally free sheaf $F$ on $X$. We begin our investigation of reflexive sheaves with the following lemma.
\begin{lemma}
On a generic complex torus $X$ of dimension $d \ge 2$ the following conditions for a coherent sheaf $G$ on $X$ are equivalent.
\begin{enumerate}
\item[(a)] $G$ is reflexive,
\item[(b)] $\Hom(T,G)=\Ext^1(T,G)=0$ for all torsion sheaves $T$.
\end{enumerate}
\end{lemma}
\begin{proof} (a) $\Longrightarrow$ (b) For any extension $0 \rightarrow G \rightarrow F \rightarrow T\rightarrow 0$ of a torsion sheaf $T$ by $G$ we consider the commutative diagram
\[ \xymatrix @C=1cm @R=1cm { 0 \ar[r] & G \ar[r]^i \ar[d]_\gamma^\wr & F \ar[r] \ar[d]^\alpha & T \ar[r] \ar[d] & 0 \\ 0 \ar[r] & G^{\vee\vee} \ar[r]^{i^{\vee\vee}} & F^{\vee\vee} \ar[r] & T'\ar[r] & 0 } \]
with exact rows and a suitable torsion sheaf $T'$. Since $G^{\vee\vee}$ and $F^{\vee\vee}$ are reflexive, the morphism $i\dd$ is determined on a complement of a Zariski-closed subset $Z$ of codimension $\ge 2$. If we take $Z=\Supp(T')$, we see that $i^{\vee\vee}:G^{\vee\vee}\longrightarrow F^{\vee\vee}$ is an isomorphism. The morphism $\pi\mathrel{\mathop{:}}=\gamma^{-1}\circ(i^{\vee\vee})^{-1}\circ\alpha$  splits our extension. The vanishing $\Hom(T,G)=0$ is obvious, because $G$ is torsionfree.\\
(b) $\Longrightarrow$ (a) Like every coherent sheaf, $G$ fits into an exact sequence
\[ 0\longrightarrow S \longrightarrow G \longrightarrow G\dd \longrightarrow T \longrightarrow 0\]
with torsion sheaves $S$ and $T$. Due to our assumption $S=0$ and the resulting short exact sequence splits. But the reflexive sheaf $G\dd$ has no torsion subsheaves, hence $T=0$ and $G$ is reflexive. \end{proof}

\begin{cor}
Assume $X$ is a generic complex torus of dimension $d\ge 3$. If $0\rightarrow F_1 \rightarrow F_2 \rightarrow F_3 \rightarrow 0$ is a short exact sequence in $\Coh(X)$ with a locally free sheaf $F_1$  and a reflexive sheaf $F_2$, then the sheaf $F_3$ is also reflexive. 
\end{cor}
\setcounter{cor1}{\value{prop}}
\begin{proof} Apply $\Hom(T,-)$ to the short exact sequence and apply (1) and Lemma 2.2. \end{proof}
Let $\omega$ be a K\"ahler class and denote as usual the slope
\[\frac{\int_X \ce_1(E)\wedge \omega^{d-1}}{\rk(E)} \]
of a torsionfree sheaf $E$ with $\mu_\omega(E)$. There is the notion of $\mu_\omega$-(semi)stability, and on a generic torus $X$ of dimension $d\ge2$ every torsionfree sheaf is semistable with slope $\mu_\omega(E)=0$.\\ 
The following important proposition is a special case of a theorem by Bando and Siu \cite{Bando-Siu}.
\begin{prop}
Let $X$ be a generic complex torus of dimension $d\ge 3$ with fixed K\"ahler metric $\omega$. Then every $\mu_\omega$-stable reflexive sheaf $F$ is a line bundle in $\Pic^0(X)$.
\end{prop}
\begin{proof} (see \cite{Bando-Siu} for more details) Bando and Siu construct a canonical Hermite--Einstein connection on the restriction of $F$ to an open set on which $F$ is locally free and whose complement consists of finitely many points. The curvature is $L^2$-integrable and satisfies the Bogomolov--L\"ubke inequality on $X$. Because of $\ce_1(F)=\ch_2(F)=0$ this connection is flat outside this finite set of points. Since the points have codimension $\ge2$, this flat connection has no local monodromy and one can extend the flat bundle to a flat bundle on $X$. Since $F$ is reflexive, it coincides with this flat bundle up to isomorphism.  The connection on the stable bundle $F$ corresponds to an irreducible representation of the abelian fundamental group of $X$. Thus, $F$ is a line bundle.\end{proof}

\begin{prop}
On a generic complex torus $X$ of dimension $d\ge 3$ every reflexive sheaf is locally free and admits a locally free filtration with quotients in $\Pic^0(X)$.
\end{prop}
\begin{proof} Since $\ce_1(F)=0$, every reflexive sheaf $F$ is $\mu_\omega$-semistable and admits a Jordan--H\"older filtration
\[ 0\subset F_0 \subset F_1 \subset \ldots F_n=F\]
with stable quotients. We may assume that $F_i$ is reflexive for all $0\le i\le n$. Due to the previous proposition $F_0\in \Pic^0(X)$. Furthermore, because of Corollary \arabic{section}.\arabic{cor1}, $F_1/F_0$ is reflexive. Hence $F_1/F_0\in \Pic^0(X)$. Since $F_0$ and $F_1/F_0$ are locally free, $F_1$ is also  locally free. Now we proceed in this way and obtain the assertion. \nolinebreak \end{proof}
\textbf{Remark 2.6.} Note that Proposition \arabic{section}.\arabic{prop} implies that there are nontrivial morphisms $L_1\rightarrow F$ and $F\rightarrow L_2$ with $L_1,L_2\in \Pic^0(X)$ for every reflexive sheaf $F$.\\
\\
In order to use these results, we will in the following assume $\dim X\ge 3$ .

\section{Some stability conditions on generic tori}

In this section we construct and characterize certain stability conditions on $\Der^b(X)$. Recall, a stability condition on $X$ consists of a bounded t-structure on $\Der^b(X)$ and an additive function on the K-group of its heart satisfying certain properties. \\
On $\Der^b(X)$ there is the standard t-structure with heart $\Coh(X)=\mathrel{\mathop{:}}\Coh_{(0)}(X)$. For the construction of other t-structures we follow the method of Happel, Reiten, and Smal{\o} using torsion pairs.
\begin{defin}
A torsion pair in an abelian category $\altA$ is a pair of full subcategories $(\altT,\altF)$ of $\altA$ with the property $\Hom_\altA(T,F)=0$ for $T\in \altT$ and $F\in \altF$. Furthermore, every object $E\in\altA$ fits into a short exact sequence
\[ 0\longrightarrow T\longrightarrow E\longrightarrow F\longrightarrow 0\]
for some objects $T\in\altT$ and $F\in \altF$.   
\end{defin}
For the abelian category $\Coh(X)$ the two subcategories $\altT\mathrel{\mathop{:}}=\{\mbox{torsion sheaves}\}$ and $\altF\mathrel{\mathop{:}}=\{\mbox{torsionfree sheaves}\}$ form a torsion pair. The following lemma illustrates the importance of this notion.
\begin{lemma}[\cite{HRS}, Proposition 2.1]
Suppose $\altA$ is the heart of a bounded t-structure on a triangulated category $\altD$ and let us denote by $H:\altD\rightarrow \altA$ the cohomology functor with respect to this t-structure. For every torsion pair $(\altT,\altF)$ in $\altA$ the full subcategory 
\[ \altA^\sharp=\{E\in \altD \mid H^i(E)=0 \mbox{ for } i\notin\{-1,0\}, H^{-1}(E)\in \altF \mbox{ and } H^0(E)\in \altT\} \]
is the heart of a bounded t-structure on $\altD$.
\end{lemma}
Using our torsion pair on $\Coh(X)$ we obtain a new t-structure on $\Der^b(X)$ whose heart $\Coh(X)^\sharp=\mathrel{\mathop{:}}\Coh_{(1)}(X)$ consists of complexes $E$ of length two with a torsion sheaf $H^0(E)$ and a torsionfree sheaf $H^{-1}(E)$. \\
We claim that on a generic torus $X$ of dimension $d\ge 3$ the pair $\altT_{(1)}=\altT=\{\mbox{torsion}$ sheaves$\}$ and $\altF_{(1)}=\{\mbox{locally free sheaves}\}[1]$ is a torsion pair in $\Coh_{(1)}(X)$. For a torsion sheaf $T$ and a locally free sheaf $F$ we have $\Ext^1(T,F)=0$ due to (\ref{eq}). Hence what remains to show is the existence of a short exact sequence as in the definition of a torsion pair. For any $E\in \Coh_{(1)}(X)$ there is a triangle
\[ T \longrightarrow H^{-1}(E)[1] \longrightarrow F[1] \longrightarrow T[1] \]
with locally free $F\mathrel{\mathop{:}}=H^{-1}(E)^{\vee\vee}$ and $T\mathrel{\mathop{:}}=F/H^{-1}(E)\in \altT_{(1)}$. We denote by $C$ the cone of the composition $T\rightarrow H^{-1}(E)[1] \rightarrow E$. From the octahedron axiom we get the triangle
\[ F[1] \longrightarrow C \longrightarrow H^0(E) \longrightarrow F[2] \]
and conclude $C=H^0(E)\oplus F[1]$, because $\Ext^1(H^0(E),F[1])= \Ext^2(H^0(E),F)=0$. If we define $K$ as the cone of the composition $E\rightarrow C \rightarrow F[1]$, we get the triangle
\[ K[-1] \longrightarrow E \longrightarrow F[1] \longrightarrow K .\]
Using the associated long exact cohomology sequence in $\Coh(X)$ and the definition of $F$ we see $K[-1]\in \altT_{(1)}$ and we are done. By definition the heart $\Coh_{(1)}(X)^\sharp=\mathrel{\mathop{:}}\Coh_{(2)}(X)$ of the new t-structure consists of objects $E$ which fit into a triangle
\[ (F[1])[1]=F[2] \longrightarrow E \longrightarrow T \]
with some torsion sheaf $T$ and some locally free sheaf $F$. For $\dim(X)=d> 3$ any such triangle splits and we get $E=T\oplus F[2]$. It is easy to check that $\altT_{(2)}=\altT=\{\mbox{torsion sheaves}\}$ and $\altF_{(2)}=\{\mbox{locally free sheaves}\}[2]$ define a torsion pair on $\Coh_{(2)}(X)$. For every object $E$ of the new abelian category $\Coh_{(2)}(X)^\sharp=\mathrel{\mathop{:}}\Coh_{(3)}(X)$ one has a triangle
\[ (F[2])[1]=F[3] \longrightarrow E\longrightarrow T \]
with a torsion sheaf $T$ and some locally free sheaf $F$. For $d>4$ we proceed in this way. Eventually one has $d$ bounded t-structures with hearts $\Coh_{(p)}(X), 0\le p<d$. In the case $0<p$ every object $E\in \Coh_{(p)}(X)$ fits into a unique triangle
\[ F[p] \longrightarrow E\longrightarrow T\]
with some torsion sheaf $T=H^0(E)$. The sheaf $F=H^{-p}(E)$ is torsionfree and, moreover, locally free for $p\ge 2$. In the case $2\le p<d-1$ the extension is trivial.    
\begin{lemma}
For every $0\le p < d$ the category $\altT$ of torsion sheaves is an abelian subcategory of $\Cohp$, i.e. if $f:S\rightarrow T$ is a morphism in $\altT$ and if we denote the kernel of $f$ in $\altT$ and in $\Cohp$ by $\ker f$ resp. $\ker_{(p)}f$, then $\ker f=\ker_{(p)} f$ and similar for the cokernels.
\end{lemma}
\setcounter{claim}{\value{prop}}
\begin{proof} Let us denote by $H^i_{(p)}$ the $i$-th cohomology functor of the t-structure corresponding to $\Cohp$. We assume $p \ge 1$ and form the triangle $ S \xrightarrow{f} T \rightarrow M \rightarrow S[1]$. Then we have $H^0(M)=\coker f$ and $H^{-1}(M)=\ker f$ as well as $H^0_{(p)}(M)=\coker_{(p)}f=\mathrel{\mathop{:}}C$ and $H^{-1}_{(p)}(M)=\ker_{(p)}f=\mathrel{\mathop{:}}K$. We form the long exact cohomology sequence in $\Coh(X)$ of the triangle $K[1] \rightarrow M \rightarrow C \rightarrow K[2]$ and use $K,C\in \Cohp$.
\begin{eqnarray*}
&& 0 \longrightarrow H^{-p}(K) \longrightarrow \underbrace{H^{-p-1}(M)}_{=0} \longrightarrow 0 \longrightarrow H^{1-p}(K) \longrightarrow \underbrace{H^{-p}(M)}_{\mbox{\scriptsize  torsion}} \longrightarrow    \\ && \underbrace{H^{-p}(C)}_{\mbox{\scriptsize torsionfree}} \longrightarrow  \underbrace{H^{2-p}(K)}_{\mbox{\scriptsize torsion}} \longrightarrow \ldots  \longrightarrow \underbrace{H^{-2}(M)}_{=0} \longrightarrow H^{-2}(C)   \longrightarrow \underbrace{H^0(K)}_{\mbox{\scriptsize torsion}} \\ && \longrightarrow \underbrace{H^{-1}(M)}_{=\ker f, \mbox{\scriptsize  torsion}} \longrightarrow H^{-1}(C) \longrightarrow 0 \longrightarrow \underbrace{H^0(M)}_{=\coker f} \longrightarrow H^0(C)  \longrightarrow 0 
\end{eqnarray*}
From this sequence we deduce $H^{-p}(K)=0$ and $H^{-p}(C)=0$. Hence $K\cong H^0(K)$ and $C\cong H^0(C)$ and, therefore, $K=\ker f$ as well as $C=\coker f$. \end{proof}
\begin{lemma}
Any morphism $f\in \Hom(F,G)$ between torsionfree sheaves $F$ and $G$ defines a  morphism $f[1]:F[1]\longrightarrow G[1]$ in $\Coh_{(1)}(X)$ and if we denote by $\Gamma(E)$ the torsion subsheaf of a sheaf $E$, we get
$H^{-1}(\ker_{(1)}(f[1]))=\ker f , \; H^0(\ker_{(1)}(f[1]))=\Gamma(\coker f),\; H^{-1}(\coker_{(1)}(f[1]))=\coker f/\Gamma (\coker f)$ 
as well as $H^0(\coker_{(1)}(f[1])) \\ =0$. 
\end{lemma}
\begin{proof} We imitate the proof of the previous lemma. Let $M$ be defined by the triangle $F[1]\xrightarrow{f} G[1] \rightarrow M \rightarrow F[2]$. Thus, $H^{-1}(M)=\coker f$ and $H^{-2}(M)=\ker f$ are the only nontrivial cohomology sheaves. The rest of the proof is left to the reader. \end{proof}
\setcounter{claim4}{\value{prop}}
\begin{prop}
For $1\le p<d$ the abelian category $\Coh_{(p)}(X)$ is of finite length, i.e. noetherian and artinian. 
\end{prop}
\begin{proof} We show that $\Cohp$ is noetherian. The proof for $\Cohp$ been artinian is similar and left to the reader. \\ 
Take an infinite sequence $E=E_0 \sur E_1 \sur E_2 \sur \ldots$ of quotients. 
We obtain the commutative diagram
\[ \xymatrix  @C=1cm @R=1cm { E_n \ar@{->>}[d] \ar@{->>}[r] & E_{n+1} \ar@{->>}[d]  \\  H^0(E_n) \ar[r] & H^0(E_{n+1}) } \]
which shows that $H^0(E_n) \rightarrow H^0(E_{n+1})$ is an epimorphism in $\Cohp$ for all $n\ge 0$. 
%The epimorphism $E_n\sur E_{n+1}$ yields an epimorphism $H^0(E_n)\sur H^0(E_{n+1})$ in $\Cohp$ for every $n$. 
Since there are only finitely many quotients of the torsion sheaf $H^0(E)$ in $\Coh(X)$ and by Lemma \arabic{section}.\arabic{claim} also in $\Cohp$, we get $H^0(E_n)\cong H^0(E_{n+1})$ for all $n\gg 0$. Then we apply the snake lemma to 
\begin{equation} \label{eq4}
\xymatrix @C=1cm @R=1cm { 0 \ar[r] & {H^{-p}(E_n)[p]} \ar[d] \ar[r] & E_n \ar@{->>}[d] \ar[r] & H^0(E_n) \ar[d]^\wr \ar[r] & 0 \\ 0 \ar[r] & {H^{-p}(E_{n+1})[p]} \ar[r] & E_{n+1} \ar[r] & H^0(E_{n+1}) \ar[r] & 0} 
\end{equation}
which yields that $H^{-p}(E_n)[p] \sur H^{-p}(E_{n+1})[p]$ is an epimorphism. Since the rank function $\rk$ is additive, the sequence $(\rk H^{-p}(E_n))_{n\in\mathbb{N}}=((-1)^p\rk (E_n))_{n\in\mathbb{N}}$ of natural numbers decreases. Thus, without loss of generality we can assume $\rk H^{-p}(E_n)=\rk H^{-p}(E_{n+1})$ for all $n\gg 0$. Hence the kernel $K_n\in \Cohp$ of the epimorphism $H^{-p}(E_n)[p] \sur H^{-p}(E_{n+1})[p]$ has rank zero and is, therefore, a torsion sheaf. In the case $2\le p<d$ there is no triangle
\[ H^{-p}(E_{n+1})[p-1] \longrightarrow K_n \longrightarrow H^{-p}(E_n)[p]\]
with $K_n\neq 0$. Hence $H^{-p}(E_n)[p] \cong H^{-p}(E_{n+1})[p]$ and (\ref{eq4}) yields $E_n\xrightarrow{\;\sim\;}E_{n+1}$ for all $n\gg 0$. \\
In the case $p=1$ set $T_n=H^{-1}(E_n)^{\vee\vee}/H^{-1}(E_n)$ and consider the commutative diagram
\begin{equation} \label{eq5} \xymatrix @C=0.8cm @R=1cm { 0 \ar[r] & T_n \ar[d]^\alpha \ar[r] & H^{-1}(E_n)[1] \ar@{->>}[d] \ar[r] & H^{-1}(E_n)^{\vee\vee}[1] \ar[d]^{\beta[1]} \ar[r] & 0 \\ 0 \ar[r] & T_{n+1} \ar[r] & H^{-1}(E_{n+1})[1] \ar[r] & H^{-1}(E_{n+1})^{\vee\vee}[1] \ar[r] & 0 } 
\end{equation}
of exact sequences in $\Coh_{(1)}(X)$. Hence $\beta[1]$ is an epimorphism in $\Coh_{(1)}(X)$ and due to Lemma \arabic{section}.\arabic{claim4} we get $\coker\beta\in\altT$. Together with $\rk H^{-1}(E_n)=\rk H^{-1}(E_{n+1})$ this shows $\ker \beta=0$. Since the sheaves are locally free, $\coker \beta=0$, because there are no divisors. Using Lemma \arabic{section}.\arabic{claim4} we conclude $\ker_{(1)}(\beta[1])=\coker_{(1)}(\beta[1])=0$ and $\beta[1]$ is an isomorphism. Hence $\alpha$ is an epimorphism in $\Coh_{(1)}(X)$ and due to Lemma \arabic{section}.\arabic{claim} also in $\Coh(X)$. Since $T_n$ has only finitely many quotients,  $T_n \xrightarrow{\;\sim\;} T_{n+1}$ for all $n\gg 0$. If we first apply the snake lemma to (\ref{eq5}) and then to (\ref{eq4}), we obtain isomorphisms $E_n\xrightarrow{\;\sim\;}E_{n+1}$ for all $n\gg 0$.  \end{proof} 
\begin{cor}
For $0\le p< d$ the additive function $Z_{(p)}(E)=-\ch_d(E)+(-1)^p\rk(E)\cdot i$ is a stability function on $\Cohp$, where $\ch_d(E)$ is (the integral over) the d-th Chern character of $E$. Furthermore, the pair $\sigma_{(p)}\mathrel{\mathop{:}}=(Z_{(p)},\Cohp)$ is a numerical locally finite stability condition on $D^b(X)$.
\end{cor}
\begin{proof} For the first part we remark that any $0\neq E\in \Cohp$ with $\rk(E)=0$ is a torsion sheaf supported on a finite set. For those sheaves $\ch_d(E)>0$. The second assertion is clear for $0<p<d$ due to the fact that $\Cohp$ is of finite length. For $p=0$ we only have to consider the case of an infinite decreasing sequence of subsheaves 
\[ \ldots \subseteq \altG_{n+1}\subseteq \altG_n \subseteq \ldots \subseteq \altG_0=\altG,\]
because $\Coh(X)$ is noetherian. For  $n\gg 0$ we have $\rk(\altG_{n+1})=\rk(\altG_n)$ and, therefore, $Z_{(0)}(\altG_{n+1})=Z_{(0)}(\altG_n) - Z_{(0)}(\altG_{n}/\altG_{n+1}) = Z_{(0)}(\altG_n)+\ch_d(T_n)$ with the torsion sheaf $T_n\mathrel{\mathop{:}}=\altG_{n}/\altG_{n+1}$. Hence the sequence of phases does not increase for $n\gg 0$. This shows that $Z_{(0)}$ satisfies the Harder--Narasimhan property on $\Coh(X)$. \\
The condition of locally finiteness is automatically fulfilled since the values of $Z_{(p)}$ form a discrete set. \end{proof}
\textbf{Remark.} After suitable modifications in the definition of $\Cohp$ all the previous statements of this section remain true for compact complex K\"ahler manifolds without nontrivial subvarieties like generic complex tori or general deformations of Hilbert schemes of K3 surfaces (see \cite{Meinhardt}). More precisely, $\Cohp$ is the abelian category of perverse sheaves with the constant perversity function $-p$. Bounded t-structures of perverse sheaves on algebraic varieties has been investigated by M. Kashiwara (\cite{Kashiwara}) and R. Bezrukavnikov (\cite{Bezrukavnikov}). \\
\\
The next proposition gives a rough classification of the objects $E$ in $\Cohp$ which are stable with respect to $\sigma_{(p)}$.
\begin{prop} In $\Cohp$ the sheaf $k(y)$ is stable of phase $1$ for any $y\in X$ and $L[p]$ is stable of phase $1/2$ for any $L\in\Pic^0(X)$. For $0<p<d-1$ these are the only stable objects in $\Cohp$. The phases of all stable objects in $\Coh(X)$ are contained in $(0,1/2]\cup\{1\}$ and the phases of all stable objects in $\Coh_{(d-1)}(X)$ are contained in $[1/2,1]$. 
\end{prop}
\setcounter{cor2}{\value{prop}}
\begin{proof} The case $p=0$: It is an easy calculation to check the stability of $L$ for any $L\in\Pic^0(X)$ and of $k(y)$ for any $y\in X$. If $E\in \Coh(X)$ is stable but not torsion, it must be torsionfree. Otherwise there is a nontrivial morphism $k(y)\rightarrow E$ which cannot exist. Furthermore, there is a nontrivial morphism $E\rightarrow E\dd \rightarrow L$ for some $L\in\Pic^0(X)$ (see Remark 2.6). Hence $\phi(E)\le \phi(L)=1/2$.  \\
\\
The case $0<p<d-1$: For $1<p<d$ we know that $H^{-p}(E)$ is locally free for any $E\in \Cohp$. This also holds for every stable object $E\in \Coh_{(1)}(E)$ which is not a torsion sheaf. Indeed, if $H^{-1}(E)$ is not locally free, there is a nonzero morphism $T \rightarrow H^{-1}(E)[1]\rightarrow E$ coming from the extension $0\rightarrow H^{-1}(E) \rightarrow H^{-1}(E)\dd \rightarrow T \rightarrow 0$ with $T\in \altT$. This contradicts the stability of $E$. Hence $H^{-p}(E)$ is locally free for any stable $E\in\Cohp, E\notin\altT$. Due to formula (\ref{eq})
\[\Ext^1(H^0(E),H^{-p}(E)[p])= \Ext^{1+p}(H^0(E),H^{-p}(E))=0 \]
and, therefore, $E\cong H^0(E)\oplus H^{-p}(E)[p]$. Hence $E\cong H^{-p}(E)[p]$ and the only stable objects are of the form $k(y)$ with phase 1 or $F[p]$ with $F$ being locally free and with phase 1/2. 
For any $L\in \Pic^0(X)$ the complex $L[p]$ has phase $1/2$. Thus, the stable factors of $L[p]$ are of the form $F[p]$ with $F$ being locally free. Since $\rk(L[p])=(-1)^p$, the complex $L[p]$ is already stable. Conversely, due to the existence of nontrivial morphisms $L[p]\rightarrow F[p]$ any stable object has rank $(-1)^{p}$ and the assertion follows.\\
\\
The case $p=d-1$: One has $Z_{(d-1)}(E)= -\ch_d(H^0(E)) + \rk(H^{1-d}(E))\cdot i$ for any $E\in \Coh_{(d-1)}(X)$. Hence $\phi(E)\in [1/2,1]$ for all $E\in \Coh_{(d-1)}(X)$. Since the phases of $k(y)$ and of $L[d-1]$ are in the boundary of the interval $[1/2,1]$ for any $y\in X$ and $L\in \Pic^0(X)$, these objects have to be semistable. They are also stable, because their Chern character is primitiv. \end{proof}
\setcounter{claim2}{\value{prop}}
Note that any ideal sheaf $\altI_{\{p_1,\ldots,p_n\}}$ is also stable in $\Coh(X)$. Hence there is no positive lower bound for the phases of stable objects in $\Coh(X)$. Similarly, there is a sequence of stable objects in $\Coh_{(d-1)}(X)$ whose phases form a strictly increasing sequence converging to 1.
\begin{cor}
For any $0<p\le d-1$ and any $\gamma\in (0,1/2)$ the pair 
\[\sigma_{(p)}^\gamma\mathrel{\mathop{:}}=\Big(Z_{(p)}^\gamma(\cdot)=-\ch_d(\cdot)-(-1)^p\cot(\pi\gamma)\rk(\cdot),\Cohp\Big)\] 
is a numerical locally finite stability condition.
\end{cor}
\begin{proof} Since $\Cohp$ is of finite type, we only have to show $Z_{(p)}^\gamma(E)<0$ for all $E\in \Cohp$. It is enough to check this for those objects in $\Cohp$ which are stable with respect to $\sigma_{(p)}$. Using Proposition \arabic{section}.\arabic{claim2} this is an easy calculation which is left to the reader. \end{proof}
Next, consider the $\Glt$-orbits through the stability conditions $\sigma_{(p)} = (Z_{(p)},\\ \Cohp)$ and $\sigma_{(p)}^\gamma=(Z_{(p)}^\gamma,\Cohp)$ in $\stab(X)$. It is an easy exercise to check that they are disjoint. \\
At the end of this section we will characterize the set
\[ U(X)\mathrel{\mathop{:}}= \bigcup_{0\le p<d}\sigma_{(p)}\cdot \Glt  \;\;\;\cup \bigcup_{1\le p<d \atop \gamma \in (0,1/2)} \sigma_{(p)}^\gamma \cdot\Glt  \]
of our stability conditions.
\begin{prop}
Assume $X$ is a generic complex torus of dimension $d\ge 3$. If $\altP$ is a slicing on $\Der^b(X)$ with the property $k(y)\in \altP(1)$ for all $y\in X$ and $L\in \altP(\psi)$ for all $L\in \Pic^0(X)$ for a fixed $\psi\in \mathbb{R}$, then $\altA\mathrel{\mathop{:}}=\altP((0,1])=\Coh_{(p)}(X)$, where $p\in \mathbb{N}$ is the unique number with $\psi+p\in(0,1]$.
\end{prop}
\setcounter{cor4}{\value{prop}}
\begin{proof} Since $\Hom(\ho_X,k(y))\neq 0$ and $\Hom(k(y),\ho_X[d])\neq 0$, we conclude
\[ \psi\in (1-d,1)\qquad\mbox{and, therefore,} \qquad 0\le p<d.\] 
The case $p=0$: In this case $k(y)\in \altA$ for all $y\in X$ and $L\in \altA$ for all $L\in \Pic^0(X)$. Furthermore, $E\in \altP([0,1))$ for all $\sigma_{(0)}$-stable torsionfree $E\in \Coh(X)$. Indeed, for such $E$ there is a triangle
\[ (E\dd/E)[-1] \longrightarrow E \longrightarrow E\dd \longrightarrow E\dd/E\]
with the locally free sheaf $E\dd\in \altP(\psi)$ and the torsion sheaf $E\dd/E\in \altP(1)$. This shows $E\in \altP([0,1))$ If $E\notin\altP((0,1))$, we find a nontrivial morphism $E\rightarrow T[-1]$ with stable $T\in \altP(1)$. We show $T\cong k(y)$ for some $y\in X$ which contradicts $\Hom(E,k(y)[-1])=0$. \\
In order to show $T\cong k(y)$, assume $H^m(T)\neq 0$ and $H^n(T)\neq 0$ but $H^k(T)=0$ $\forall \:k\notin[m,n]$ for two integers $m\le n$. If $H^m(T)$ is not torsionfree, there are  nontrivial compositions
\[ k(y)[-m]\longrightarrow H^m(T)[-m]\longrightarrow T \quad\mbox{and}\quad T\longrightarrow H^n(T)[-n] \longrightarrow k(z)[-n]\]
for  suitable $y,z\in X$. If $H^m(T)$ is torsionfree but not reflexive, we replace the first composition by
\[ k(y)[-1-m]\longrightarrow H^m(T)[-m]\longrightarrow T \]
and if $H^m(T)$ is reflexive, we take
\[ L[-m]\longrightarrow H^m(T)[-m] \longrightarrow T\]
for a suitable $L\in \Pic^0(X)$. If $T\in \altP(1)$ is not isomorphic to $k(y)$, we get in all cases $-m\le\psi-m <1<1-n$ and, therefore, $n< 0 \le m$, a contradiction to $m\le n$. \\ 
Thus, any $\sigma_{(0)}$-stable sheaf is contained in $\altA=\altP((0,1])$ and we get $\Coh(X)\subseteq \altA$. By standard arguments $\Coh(X)=\altA$. \\
\\
The case $0<p< d$: From the proof of Proposition \arabic{section}.\arabic{cor2} we know that any $\sigma_{(p)}$-stable object $E\in \Cohp$ fits into a triangle
\[ H^{-p}(E)[p] \longrightarrow E \longrightarrow H^0(E) \]
with locally free $H^{-p}(E)\in \altP(\psi)$ and the torsion sheaf $H^0(E)\in \altP(1)$. Since $\altP(\psi+p)\subseteq \altA$ and $\altP(1)\subseteq \altA$, we see $E\in \altA$ and, therefore, $\Cohp\subseteq \altA$. Again we can conclude $\Cohp=\altA$.\end{proof}
\begin{prop}
Assume $(Z',\Cohp)$ is a locally finite numerical stability condition on $X$ with $0\le p\le d-1$ and $\phi'(k(y))=1$ for all $y\in X$. Then there is a matrix $G\in \Gl2$ with $G\cdot Z_{(p)}=Z'$ or $G\cdot Z_{(p)}^\gamma=Z'$ for a unique $\gamma\in (0,1/2)$.
\end{prop}
\setcounter{cor3}{\value{prop}}
\begin{proof} Since $Z'$ is numerical, we get $Z'(E)=-e\ch_d(E) - (-1)^pf\rk(E) + (g\ch_d(E)+(-1)^ph\rk(E))\cdot i$ for a suitable matrix 
\[ \Z{e}{f}{g}{h}\in \Mat2.\]
Since $\phi'(k(y))=1$, we obtain $g=0$ and $e>0$. If $Z'$ takes values in $(-\infty,0)$, then $h=0, f>0$ and 
\[ \Vek{\Rea Z'}{\Ima Z'}=\Z{e}{0}{0}{1}\cdot\Vek{-\ch_d-(-1)^p\cot(\pi\gamma)\rk}{0}\quad\mbox{with } \cot(\pi\gamma)=f/e.\]
This can only occur for $0<p\le d-1$ since $\Coh(X)$ is not of finite type. If the image of $Z'$ is not contained in $(-\infty,0)$,then $h>0$ and 
\[ \Vek{\Rea Z'}{\Ima Z'}=\Z{e}{-f}{0}{h}\cdot\Vek{-\ch_d}{(-1)^p\rk}.  \]
\end{proof}
Using these two propositions we get the main result of this section which characterizes the set $U(X)$ of stability conditions.
\begin{theorem}
Assume $X$ is a generic complex torus of dimension $d$. The set 
\[ U(X)= \bigcup_{0\le p<d}\sigma_{(p)}\cdot \Glt \;\;\;\cup  \bigcup_{1\le p<d \atop \gamma \in (0,1/2)} \sigma_{(p)}^\gamma \cdot\Glt  \]
is the set of all numerical locally finite stability conditions $\sigma=(Z,\altP)$ such that there exist certain real numbers $\phi$ and $\psi$ such that $k(y)\in \altP(\phi)$ for all $y\in X$ and $L\in \altP(\psi)$ for all $L\in\Pic^0(X)$.
\end{theorem}
\begin{proof} Choose some stability condition $\sigma=(Z,\altP)$ with the property described in the theorem. After applying some $G\in \Glt$ we can assume $k(y)\in \altP(1)$ $\forall \:y\in X$. Using Proposition \arabic{section}.\arabic{cor4} and Proposition \arabic{section}.\arabic{cor3} we get $\sigma\in U(X)$. Of course, every stability condition in $U(X)$ has the charaterizing property. \end{proof}

\section{The topology of $U(X)$}

In this section we study the topology of $U(X)$. As we will see, $U(X)$ is a simply connected component of $\stab(X)$.\\ The first part of this section is a more general consideration of $\Glt$-orbits in the space $\stab(\altD)$ of locally finite stability conditions on a triangulated category $\altD$. In the second part we come back to the case $\altD=\Der^b(X)$. \\ 
\\
Let $\Sigma\subseteq \stab(\altD)$ be a connected component and let us denote by $V(\Sigma)$ the linear subspace in $\Hom(K(\altD),\mathbb{C})$ such that the forgetting map 
\[\altZ:\stab(\altD)\supseteq \Sigma \ni \sigma=(Z,\altP)\longmapsto Z\in V(\Sigma)\subseteq \Hom(K(\altD),\mathbb{C})\] 
is a local homeomorphism. Given a stability condition $\sigma=(Z,\altP)\in \Sigma$ the space $V(\Sigma)$ is characterized by 
\[ V(\Sigma)=\{U\in \Hom(K(\altD),\mathbb{C})\mid  \|U\|_\sigma < \infty\}, \] 
where
\[ \| U \|_\sigma \mathrel{\mathop{:}}= \sup\left\{ \left. \frac{|U(E)|}{|Z(E)|} \;\right| E \mbox{ semistable in } \sigma\right\} \]
and $\| \cdot\|_\sigma$ can be used to define the topology on $V(\Sigma)$ \cite{Bridgeland02}. It follows that the evaluation map $V(\Sigma)\ni U \mapsto U(E)\in \mathbb{C}$ is continuous for a fixed $E\in \altD$.\\
The universal cover $\Glt$ of $\Gl2$ acts on $\stab(\altD)$ from the left by $g\cdot \sigma \mathrel{\mathop{:}}= \sigma \cdot g^{-1}$, where the the latter action is the one considered by Bridgeland in \cite{Bridgeland02}. Furthermore, there is an action from the left of the ring $\Mat2$ on $\Hom(K(\altD), \mathbb{C})$ and the map $\altZ$ commutes with these actions. Let us consider a stability condition $\sigma=(Z,\altP)\in \Sigma$ such that the image of the central charge is not contained in a real line in $\mathbb{C}$ and $\altP(1)\neq \{0\}$. We are interested in the boundary points of the orbit
\[ \sigma\cdot \Glt =(\Glt)^ {-1}\cdot\sigma =\Glt\cdot\,\sigma \subseteq \Sigma. \]
This orbit is a real submanifold of $\Sigma$ of real dimension four. It follows from the definition that the central charges of all stability conditions of this orbit factorize over the quotient by $K(\altD)_{\mathbb{R},\sigma}^\perp\mathrel{\mathop{:}}=\{ e\in K(\altD)_\mathbb{R}  \mid Z(e)=0 \}$ of real codimension two, i.e. they are contained in the closed real four-dimensional subspace $V(\Sigma)_\sigma\mathrel{\mathop{:}}=\Hom_\mathbb{R}(K(\altD)_\mathbb{R}/K(\altD)_{\mathbb{R},\sigma}^\perp,\mathbb{C})$. The map $\Mat2\ni M \mapsto M \circ Z \in V(\Sigma)$ is an $\mathbb{R}$-linear isomorphism onto $V(\Sigma)_\sigma$. This isomorphism identifies $\Gl2$ with $\altZ(\Glt\cdot\sigma)$. We write $Z(E)=\Re(E)+i\cdot\Im(E)$ with linear independent $\Re$ and $\Im\in \Hom(K(\altD),\mathbb{R})$. \\
Let us denote by $\sigma'=(Z',\altP')$ a boundary point of the orbit $\sigma\cdot\Glt$. Since the evaluation map is continuous, $Z'$ still factorizes over $K(\altD)_{\mathbb{R},\sigma}^\perp$. After applying some element of $\Glt$ to $\sigma'$, we can, therefore, assume $Z'=\Re -\cot(\pi\gamma) \Im$ with a suitable $\gamma\in (0,1)$, because semistability is a closed property and, therefore, $Z'(E)\neq 0 \;\forall E\in \altP(1)$. The line $Z'=0$ in $\mathbb{C}\cong \mathbb{R}\Re \oplus \mathbb{R}\Im$ is given by the equation $\Re = \cot(\pi \gamma)\Im$ and since $Z'(E)\neq 0\;\forall \:E$ semistable in $\sigma$, we have 
\begin{equation} \label{gammabed1} \gamma\neq \phi(E) \quad \forall \;E\;\mbox{ stable in } \;\sigma. 
\end{equation}
The following result was already known to the experts (see for example \cite{Bridgeland03} and \cite{Bridgeland06}).
\begin{prop}
The heart $\altP'((0,1])=\altP'(1)$ of $\sigma'$ is the tilt $\altA^\sigma_\gamma$ of $\altA\mathrel{\mathop{:}}=\altP((0,1])$ with respect to the torsion theory $(\altP((\gamma,1]),\altP((0,\gamma)))$, i.e.
\[ \altP'(1)=\{ E\in \altD\mid H^0(E)\in \altP((\gamma,1]), H^{-1}(E)\in \altP((0,\gamma)), H^k(E)=0\;\mbox{else }\}, \]
where $H$ denotes the cohomology functor associated to the bounded t-structure with heart $\altA$.
\end{prop}
\begin{proof} Due to (\ref{gammabed1}) the pair $(\altP((\gamma,1]),\altP((0,\gamma)))$ is indeed a torsion theory in $\altA=\altP((0,1])$ and since $Z'(E)\neq 0\;\forall \:E$ semistable in $\sigma$, we obtain $E\in \altP'(0)\;\forall \;E$ semistable in $\sigma$ with $\phi(E)\in (0,\gamma)$. Therefore, $\altP((0,\gamma))\subseteq \altP'(0)$ and, similarly, $\altP((\gamma,1])\subseteq \altP'(1)$. Hence $\altP'(1)$ contains the tilt of $\altP((0,1])$ with respect to the above torsion theory. By standard arguments one concludes equality. \end{proof}
In order to show the nonexistence of boundary points $\sigma'$,  we introduce the following two phases for our stability condition $\sigma=(Z,\altP)$ and the real number $\gamma\in(0,1)$.
\begin{eqnarray*}
&&\gamma^+\mathrel{\mathop{:}}=\inf\{\phi(E)\mid E \in \altD\mbox{ stable in }\sigma, \phi(E) > \gamma\}, \\
&&\gamma^-\mathrel{\mathop{:}}=\sup\{\phi(E)\mid E\in \altD \mbox{ stable in }\sigma, \phi(E) < \gamma\}. 
\end{eqnarray*}

Clearly $\gamma^-\le \gamma\le \gamma^+$ and there is no $E\in \altD$, stable in $\sigma$, with $\phi(E)\in (\gamma^-,\gamma^+)$. Hence for all $\gamma'\in[\gamma^-,\gamma^+]$ satisfying (\ref{gammabed1}) we obtain $\gamma^+=\gamma'^+$ and $\gamma^-=\gamma'^-$. Note that for $\gamma\in (\gamma^-,\gamma^+)$ the condition (\ref{gammabed1}) is always fulfilled.
\begin{prop}
If $\altP(\gamma^+)=\{0\}$ or $\altP(\gamma^-)=\{0\}$ and $\gamma'\in [\gamma^-,\gamma^+]$ satisfying (\ref{gammabed1}), there is no boundary point of $\sigma\cdot\Glt$ with central charge $Z'\mathrel{\mathop{:}}=Z_{\gamma'}^\sigma\mathrel{\mathop{:}}=\Re -\cot(\pi\gamma')\Im$.
\end{prop}
\setcounter{claim2}{\value{prop}}
\begin{proof} We consider the case $\altP(\gamma^+)=\{0\}$. The second case is similar. If there is a boundary point $\sigma'=(Z'=Z_{\gamma'}^\sigma,\altP')$, we can replace $Z_{\gamma'}^\sigma$ and assume $\gamma'=\gamma'^+=\gamma^+$. Indeed, since $\altP'(1)$ is of finite length and $\gamma^+$ satisfies (\ref{gammabed1}), the pair $\sigma^+\mathrel{\mathop{:}}=(Z_{\gamma^+}^\sigma,\altP')$ is a locally finite stability condition. It is easy to see that $\sigma^+$ is still in the boundary of $\sigma\cdot \Glt$. Since $\altZ$ is a local homeomorphism, there is an open neighbourhood of $\sigma^+$ in $\Sigma$ which is isomorphic to an open ball in $V(\Sigma)$. The intersection of this ball with $V(\Sigma)_\sigma$ can be identified with an open ball in $\Mat2$ with center
\[ Z_{\gamma^+}^\sigma  \cong \Z{1}{-\cot(\pi\gamma^+)}{0}{0}. \]
Such a ball contains the central charge $Z_{\gamma''}^\sigma=\Re -\cot(\pi\gamma'')\Im$ with $\gamma''\in (\gamma^+,\gamma^+ + \varepsilon)$ and $\varepsilon > 0$ sufficiently small. By definition of $\gamma^+$ we can assume without loss of generality $Z_{\gamma''}^\sigma(E)=0$ for some $\sigma$-stable $0\neq E\in \altD$ . As $Z_{\gamma''}^\sigma$ is a boundary point of the orbit $\Gl2\cdot Z=(\Gl2)^{-1}\cdot Z=\altZ(\sigma\cdot\Glt)$ and semistability is a closed property, this $E$ is still semistable in the stability condition lying in the neighbourhood of $\sigma^+$ and mapped by $\altZ$ onto $Z_{\gamma''}^\sigma$. This contradicts $Z_{\gamma''}^\sigma(E)=0$. \end{proof}
Due to Proposition \arabic{section}.\arabic{prop}, we have to assume $\altP(\gamma^+)\neq\{0\}$ and $\altP(\gamma^-)\neq\{0\}$ in order to obtain stability conditions with central charges $Z_\gamma^\sigma$ in the boundary of the orbit $\sigma\cdot\Glt$. \\
\\
\\
\begin{center}
\setlength{\unitlength}{1.0cm}
\begin{picture}(6.0,6.0)
\put(-0.3,3){\vector(1,0){6.6}}
\put(6.0,2.6){\begin{math}\Re \end{math}}
\put(3,-0.3){\vector(0,1){6.6}}
\put(2.6,6.1){\begin{math}\Im \end{math}}
\multiput(0,3)(1,0){3}{\circle*{0.1}}
\multiput(3.5,3.2)(0.8,0){4}{\circle*{0.1}}
\multiput(0,3.8)(0.5,0){4}{\circle*{0.1}}
\multiput(3.2,4.1)(1,0){4}{\circle*{0.1}}
\multiput(0.4,4.8)(0.8,0){4}{\circle*{0.1}}
\multiput(1.2,5.4)(1.3,0){5}{\circle*{0.1}}
\multiput(4.5,4.8)(1,0){3}{\circle*{0.1}}
\multiput(0,6)(1.1,0){4}{\circle*{0.1}}
\multiput(5.4,6)(1.1,0){2}{\circle*{0.1}}
\put(2,0){\line(1,3){2.1}}
\put(0.6,0){\line(4,5){5.1}}
\put(1.0,0){\line(2,3){4.3}}
\put(0,2.2){\begin{math} \altZ_{\gamma'}^\sigma<0 \end{math}}
\put(3.2,6.3){\begin{math} \altZ_{\gamma^+}^\sigma=0 \end{math}}
\put(-0.7,0.1){\begin{math} \altZ_{\gamma^-}^\sigma=0 \end{math}}
\put(4.6,6.5){\begin{math} \altZ_{\gamma}^\sigma=0 \end{math}}
\put(4.5,2,2){\begin{math} \altZ_{\gamma'}^\sigma>0 \end{math}}
\end{picture} \\
\vspace{0.5cm}
\scriptsize The dots are the central charges of the $\sigma$-semistable objects in $\altA$.\\
\vspace{0.4cm}
\end{center}
As in the end we want to avoid boundary points, we need a criterion that excludes the cases $\altP(\gamma^+)\neq \{0\}$ and $\altP(\gamma^-)\neq \{0\}$. This is only possible in special situations and the following will be enough in the geometric context we are interested in.
\begin{lemma}
Suppose there exists a sequence $E_n\in \altP(\gamma^+), n\in \mathbb{N}$, of non isomorphic simple objects. Then there is no object $I\in \altP((0,\gamma^-])$ with $\Ext^1(E_n,I)\neq 0$ for all $n\in \mathbb{N}$.
\end{lemma}
\setcounter{claim3}{\value{prop}}
\begin{proof}
If such an object $I$ exists, we construct by induction a sequence of nontrivial extensions
\[ 0 \longrightarrow I_n \longrightarrow I_{n+1}\longrightarrow E_n \longrightarrow 0\]
in $\altA=\altP((0,1])$ with $I_n\in \altP((0,\gamma^-])$ and the additional property $\Ext^1(E_k,I_n)\neq0$ for all $k\ge n$ and $n\in \mathbb{N}$. Since $Z(I_{n+1})=Z(I_n)+Z(E_n)$, we get $\phi(I_n)>\gamma^-$ for $n\gg 0$ which contradicts $I_n\in \altP((0,\gamma^-])$. \\
The construction of $I_n$ starts with $I_0=I$. Due to our assumption this is possible. Assume we have constructed $I_n\in \altP((0,\gamma^-])$. Choose an element $0\neq e\in \Ext^1(E_n,I_n)$ and consider the corresponding nontrivial extension in $\altA$
\[ 0 \longrightarrow I_n \longrightarrow I_{n+1}\longrightarrow E_n \longrightarrow 0.\]
For any $0\neq F\in P(\gamma^-,1]=P[\gamma^+,1]$ stable in $\sigma$ we get the following long exact sequence
\[ 0 \longrightarrow \Hom(F,I_{n+1}) \longrightarrow \Hom(F,E_n) \longrightarrow \Ext^1(F,I_n) \longrightarrow \Ext^1(F,I_{n+1}). \]    
Now, $\Hom(F,E_n)=0$ unless $F=E_n$ and in the latter case $\Hom(E_n,E_n)= \mathbb{C}\cdot\Id_{E_n}$. But $\Id_{E_n}$ is mapped to $0\neq e\in \Ext^1(E_n,I_n)$. Therefore, $\Hom(F,I_{n+1}) =0$ for all $F\in \altP((\gamma^-,1])$ and we conclude $I_{n+1}\in \altP((0,\gamma^-])$. Furthermore, the map $\Ext^1(E_k,I_n) \longrightarrow \Ext^1(E_k,I_{n+1})$ is an injection for $k\ge n+1$. Hence $\Ext^1(E_k,I_{n+1}) \\ \neq 0$ for all $k\ge n+1$ by the induction hypothesis and we are done. 
\end{proof}
Using this we get our main result of this section.
\begin{theorem}
Assume $X$ is a generic complex torus of dimension $d\ge 3$. Then 
\[ U(X)\mathrel{\mathop{:}}=\bigcup_{0\le p<d}\sigma_{(p)}\cdot \Glt \;\;\;\cup  \bigcup_{1\le p<d \atop \gamma \in (0,1/2)} \sigma_{(p)}^\gamma \cdot\Glt  \]
is a simply connected component of $\stab(X)$ but $\altZ:U(X) \longrightarrow \altZ(U(X))$ is not a covering. 
\end{theorem}
\begin{proof} On a generic complex torus of dimension $d$ any stability function of a numerical stability condition is a complex linear combination of $\ch_0=\rk$ and $\ch_d$. Since the orbits $\sigma_{(p)}\cdot\Glt$ are of real dimension four, they are open in $\stab(X)$. We describe the closure of these open orbits beginning with that of $\sigma_{(0)}=(Z_{(0)},\Coh(X))$. \\
We want to exclude boundary points with $\gamma\in (0,1/2]$. In order to apply Proposition 4.\arabic{claim2}, we show $\altP(\gamma^+)=\{0\}$. Indeed, if $0\neq E\in \altP(\gamma^+)$ is a stable sheaf, then it is torsionfree, because  $\gamma^+\le 1/2$. Now, choose a sequence of numerical trivial line bundles $L_n\in \Pic^0(X)$ with $L_m^{\rk(E)}\neq L_n^{\rk(E)}$ for all $m\neq n$. Hence $E\otimes L_m \neq E\otimes L_n$ for $m\neq n$, because of $\det(E\dd \otimes L)=\det(E\dd)\otimes L^{\rk(E)}$ for every $L\in \Pic(X)$. Furthermore, the sheaves $E_n\mathrel{\mathop{:}}=E\otimes L_n$ are also $\sigma_{(0)}$-stable of phase $\gamma^+$. We introduce the sheaf $P\mathrel{\mathop{:}}=k(y)$ for some $y\in X$. Choose an epimorphism $f:E_0 \sur P$ and denote the kernel by $I$. We prove $I\in \altP((0,\gamma^-])$ and $\Ext^1(E_n,I)\neq 0$ for all $n\in \mathbb{N}$ which contradicts Lemma 4.\arabic{claim3}. Thus, $\altP(\gamma^+)=\{0\}$. \\
In order to show $I\in \altP((0,\gamma^-])$, we take a $\sigma_{(0)}$-stable sheaf $F\in \altP((\gamma^-,1])=\altP([\gamma^+,1])$ and consider the long exact sequence
\[ 0\longrightarrow \Hom(F,I) \longrightarrow \Hom(F,E_0) \longrightarrow \Hom(F,P) \longrightarrow \Ext^1(F,I).\]
Now, $\Hom(F,E_0)=0$ unless $F=E_0$ and in the latter case $\Hom(E_0,E_0)=  \mathbb{C}\cdot\Id_{E_0}$. But $\Id_{E_0}$ is mapped to $0\neq f\in \Hom(E_0,P)$. Therefore, $\Hom(F,I)=0$ for all $F\in\altP((\gamma^-,1])$ and $I\in \altP((0,\gamma^-])$ follows. \\
For the second property of $I$ we consider the inclusion $\Hom(E_n,P)\hookrightarrow \Ext^1(E_n,I)$  and note that the former set contains $f\otimes id_{L_n} \neq 0$ for all $n\in \mathbb{N}$.\\
On the other hand, for every $\gamma\in (1/2,1)$ we obtain $\sigma_{(1)}^{1-\gamma}$ as a boundary point. \\
In the case $0<p<d-1$ the situation is very easy. There are two regions of boundary points of the orbit $\sigma_{(p)}\cdot\Glt$. For $\gamma\in (0,1/2)$ the boundary points are given by $\sigma_{(p)}^\gamma$ and for $\gamma\in (1/2,1)$ the boundary points are $\sigma_{(p+1)}^{1-\gamma}$. \\
The case $p=d-1$ is similar to the case $p=0$. First of all $E\otimes L\in \Coh_{(d-1)}(X)$ for all $E\in \Coh_{(d-1)}(X)$ and $L\in \Pic^0(X)$. Indeed, this is true for $E\cong k(y)$ and $E\cong H^{1-d}(E)[d-1]$ locally free. But any $E\in \Coh_{(d-1)}(X)$ is an extension of such special objects and tensoring with $L$ maps extensions to exensions. Furthermore, $E\otimes L\not\cong E$ for all $E\in \Coh_{(d-1)}(X)\setminus \altT$ and all $L\in \Pic^0(X)$ with $L^{\rk(E)}\not\cong \ho_X$, because $H^{1-d}(E\otimes L)=H^{1-d}(E)\otimes L\not\cong H^{1-d}(E)$.\\
Now, we can exclude boundary points with $\gamma\in (1,1/2)$ in the same way as for $\sigma_{(0)}$. Note that $\gamma^+<1$, because there are $\sigma_{(d-1)}$-stable objects with  phases sufficiently close to 1. The definitions of $P, f$ and $I$ are given by the short exact sequence
\[ 0 \longrightarrow \underbrace{H^{1-d}(E_0)[d-1]}_{\mathrel{\mathop{:}}=I} \longrightarrow E_0 \xrightarrow{=\mathrel{\mathop{:}}f} \underbrace{H^0(E_0)}_{\mathrel{\mathop{:}}=P}\longrightarrow 0.\]
Since $H^{1-d}(E_0)[d-1]$ has phase $1/2$ (see Proposition 3.\arabic{cor2}) and $\gamma^-\ge 1/2$, the property $I\in \altP((0,\gamma^-])$ is obvious in this case.\\ 
On the other hand, for every $\gamma\in (0,1/2)$ we obtain $\sigma_{(d-1)}^{\gamma}$ as a boundary point. \\   
Hence the open four-dimensional orbits $\sigma_{(p)}\cdot\Glt$ are successive connected by the three dimensional `walls' $\cup_{\gamma\in (0,1/2)} \sigma_{(p)}^\gamma\cdot\Glt$. Furthermore, the connected set $U(X)$ is closed in $\stab(X)$. The image of $U(X)$ under $\altZ:\stab(X) \longrightarrow (H^0(X,\mathbb{C})\oplus H^d(X,\mathbb{C}))^\vee\cong \Mat2$ is an open subset. Since $\altZ$ is a local homeomorphism, $U(X)$ is also open and hence a connected component of $\stab(X)$. \\ 
It is easy to see that $(\{\sigma_{(p-1)},\sigma_{(p)}\} \cup \{\sigma_{(p)}^\gamma\mid\gamma\in (0,1/2)\})\cdot\Glt$ is the universal cover of its image for all $0<p\le d-1$. This image has fundamental group $\mathbb{Z}$ which is `resolved' by the shift functor $[2]$. Using the Seifert--van Kampen theorem one concludes $\pi_1(U(X))=0$. Since the number of preimages of a stability function is not constant, $U(X)$ is not the universal cover of its image. \nolinebreak[4] \end{proof} 

\bibliographystyle{plain}
\bibliography{Literatur}
\vfill
\textsc{Mathematisches Institut, Universit\"at Bonn, Germany}\\
\textit{E-mail address:} \texttt{sven@math.uni-bonn.de}

\end{document}